\documentclass[12pt]{article}

\usepackage{amsthm, amsmath, amsfonts, amssymb}
\usepackage[dvips, all, knot]{xy}
\xyoption{all}

\usepackage{geometry}
\usepackage{fancyhdr}
\usepackage{sectsty}
\allsectionsfont{\normalsize\centering}

\newcommand{\Card}{{\rm Card}}
\newcommand{\diag}{{\rm diag}}
\newcommand{\dist}{{\rm d}}
\newcommand{\GL}{{\rm GL}}
\newcommand{\GSp}{{\rm GSp}}
\newcommand{\ord}{{\rm ord}}
\newcommand{\PGSp}{{\rm PGSp}}
\newcommand{\SL}{{\rm SL}}
\newcommand{\Sp}{{\rm Sp}}
\newcommand{\Stab}{{\rm Stab}}
\newcommand{\ti}{{\rm totally isotropic}}
\newcommand{\vol}{{\rm vol}}

\newcommand{\cB}{{\cal B}}
\newcommand{\C}{{\mathbb C}}
\newcommand{\cE}{{\cal E}}
\newcommand{\F}{{\mathbb F}}
\newcommand{\cO}{{\cal O}}
\newcommand{\Z}{{\mathbb Z}}
\newcommand{\ue}{\underline{\varepsilon}}

\newtheorem{nthm}{Theorem}[section]
\newtheorem{ncor}[nthm]{Corollary}
\newtheorem{nlem}[nthm]{Lemma}
\newtheorem{nprop}[nthm]{Proposition}

\begin{document}
\title{Expanders and the Affine Building of $\Sp_n$}
\author{Alison Setyadi}
\date{}
\maketitle
\thispagestyle{empty}

\begin{abstract}
For $n \geq 2$ and a local field $K$, let $\Delta_n$ denote the affine building
naturally associated to the symplectic group $\Sp_n(K)$. We compute the
spectral radius of the subgraph $Y_n$ of $\Delta_n$ induced by the special
vertices in $\Delta_n$, from which it follows that $Y_n$ is an analogue of a
family of expanders and is non-amenable.
\end{abstract}

\section*{Introduction}
In \cite{bien}, Bien considers the problem of constructing a model for an
efficient communication network. Such a model can be represented by a
magnifier---a graph with a small number of edges such that every subset of
vertices has many distinct neighbors (see \cite[p.\ 6]{bien}). For a finite
graph $X$, the ``quality'' of $X$ as a network can be quantified by its
isoperimetric constant $h(X)$ (see \cite[p.\ 1]{dsv}). Davidoff, Sarnak, and
Valette explore the problem of explicitly constructing a family of expanders;
i.e., a family $\{X_m\}$ of finite, connected, $r$-regular ($r \geq 2$) graphs
with $\Card(X_m) \rightarrow \infty$ as $m \rightarrow \infty$ such that there
is an $\varepsilon > 0$ satisfying $h(X_m) \geq \varepsilon$ for all $m$
\cite[Definition 0.3.3]{dsv}. Families of expanders have become building blocks
in many engineering applications, including network designs, complexity theory,
coding theory, and cryptography (see \cite[p.\ 3]{dsv} and the references cited
there). Note that by \cite[p.\ 4]{dsv}, an infinite family of $r$-regular
Ramanujan graphs is not only a family of expanders but is also optimal from the
spectral viewpoint.

Let $K$ be a local field, and let $\Xi_n$ denote the affine building naturally
associated to $\SL_n(K)$. In \cite[Example 3]{sw}, Saloff-Coste and Woess
explicitly calculate the spectral radius of the simple random walk on the
one-complex $X_n$ of $\Xi_n$. This gives the spectral radius $\rho(X_n)$ of
$X_n$ \cite[Theorem 1]{csz}. Then $h(X_n) > 0$ (by \cite[Theorem 3.3]{bms}) and
the number of vertices in $X_n$ contained in concentric balls grows
exponentially with the radius of the balls (by \cite[Theorem 2.2]{bms});
alternatively, $X_n$ is \emph{expanding}. Note that $h(X_n) > 0$ implies $X_n$
non-amenable (see the paragraph following Theorem \ref{specrad}).

Let $\Delta_n$ denote the affine building naturally associated to the
symplectic group $\Sp_n(K) \leq \SL_{2n}(K)$. In this paper, we consider the
subgraph $Y_n$ of $\Delta_n$ induced by the special vertices in $\Delta_n$ (all
the vertices in $\Xi_n$ are special). Using the techniques of \cite[Example
3]{sw}, we compute $\rho(Y_n)$. Then $h(Y_n) > 0$ (by \cite[Theorem 3.3]{bms})
and the number of vertices in $Y_n$ contained in concentric balls grows
exponentially with the radius of the balls (by \cite[Theorem 2.2]{bms}); hence,
$Y_n$ is also expanding and non-amenable.

After completing this work, we learned that it is also possible to derive the
formula in Proposition \ref{Spnspecrad} using the techniques and results of
James Parkinson. Parkinson's approach is quite general, as it takes a
building-theoretic perspective rather than the group-theoretic one we use. As
in \cite[Example 3]{sw}, Parkinson's approach is through a simple random walk:
by \cite[Theorem 6.3]{parkinsonmz} and general facts about $C^\star$-algebras,
the spectral radius of an isotropic random walk (of which a simple random walk
is an example) on an arbitrary thick, regular, affine building of irreducible
type is $\widehat{A}(1)$, where $A$ is the transition operator of the random
walk and $\widehat{A}$ its Gelfand transform. To express $\widehat{A}(1)$ for
the graph $Y_n$ in terms of the order $q$ of the residue field of $K$ (as in
Proposition \ref{Spnspecrad} below), one identifies and uses the underlying
root system of the building $\Delta_n$, together with results about the
Macdonald spherical functions defined in \cite[p.\ 580]{parkinsonmz}. In
contrast, our approach is through the natural association of $\Delta_n$ with
$\Sp_n(K)$---in particular, the transitive action of $\GSp_n(K)$ (the analogue
of $\GL_{2n}(K)$ for $\Sp_n(K)$) on the special vertices in $\Delta_n$ (see
Proposition \ref{GSpaction}). As a result, we deduce properties about
$\Sp_n(K)$ and $\GSp_n(K)$---for example, we produce a solvable subgroup of
$\GSp_n(K)$ that acts transitively on the vertices in $Y_n$. In addition, we
characterize the set of vertices in $Y_n$ adjacent to a given one in terms of
orbits (Proposition \ref{oneighbors}).

I thank Paul Garrett for his help with the Haar measure on $Q$ described in
Section \ref{comp}. Finally, the results contained here form part of my
doctoral thesis, which I wrote under the guidance of Thomas R.\ Shemanske.

\section{The Affine Building $\Delta_n$ of $\Sp_n(K)$} \label{prelims}

Fix a local field $K$ with discrete valuation ``$\ord$,'' valuation ring $\cO$,
uniformizer $\pi$, and residue field $k \cong \F_q$. The affine building
$\Delta_n$ naturally associated to $\Sp_n(K)$ can be modeled as an
$n$-dimensional simplicial complex as follows (see \cite[pp.\  336 --
337]{garrett}). Fix a $2n$-dimensional $K$-vector space $V$ endowed with a
non-degenerate, alternating bilinear form $\langle \cdot, \cdot \rangle$, and
recall that a subspace $U$ of $V$ is \emph{totally isotropic} if $\langle u, u'
\rangle = 0$ for all $u, u' \in U$. A \emph{lattice} in $V$ is a free, rank
$2n$, $\cO$-submodule of $V$, and two lattices $L$ and $L'$ in $V$ are
\emph{homothetic} if $L' = \alpha L$ for some $\alpha \in K^\times$; write $[L]$
for the homothety class of the lattice $L$. A lattice $L$ is \emph{primitive} if
$\langle L, L \rangle \subseteq \cO$ and $\langle \cdot, \cdot \rangle$ induces
a non-degenerate, alternating $k$-bilinear form on $L/\pi L$. Then a
\emph{vertex} in $\Delta_n$ is a homothety class of lattices in $V$ with a
representative $L$ such that there is a primitive lattice $L_0$ with $\langle
L, L \rangle \subseteq \pi\cO$ and $\pi L_0 \subseteq L \subseteq L_0$;
equivalently, $L/\pi L_0$ is a \ti\ $k$-subspace of $L_0/\pi L_0$. Two vertices
$t, t' \in \Delta_n$ are \emph{incident} if there are representatives $L \in t$
and $L' \in t'$ such that there is a primitive lattice $L_0$ with $\langle L, L
\rangle \subseteq \pi\cO$, $\langle L', L' \rangle \subseteq \pi\cO$, and
either $\pi L_0 \subseteq L \subseteq L' \subseteq L_0$ or $\pi L_0 \subseteq
L' \subseteq L \subseteq L_0$. Thus, a maximal simplex or \emph{chamber} in
$\Delta_n$ has $n + 1$ vertices $t_0, \ldots, t_n$ with representatives $L_i
\in t_i$ such that $L_0$ is primitive, $\langle L_i, L_i \rangle \subseteq
\pi\cO$ for all $1 \leq i \leq n$, and $\pi L_0 \subsetneq L_1 \subsetneq
\cdots \subsetneq L_n \subsetneq L_0$.

Recall that a basis $\{u_1, \ldots, u_n, w_1, \ldots, w_n\}$ for $V$ is
\emph{symplectic} if $\langle u_i, w_j \rangle$ $= \delta_{ij}$ (Kronecker
delta) and $\langle u_i, u_j \rangle = 0 = \langle w_i, w_j \rangle$ for all
$i, j$. If a $2$-dimensional, \ti\ subspace $U$ of $V$ is a hyperbolic plane,
then a frame is an unordered $n$-tuple $\{\lambda_1^1, \lambda_1^2\}, \ldots,
\{\lambda_n^1, \lambda_n^2\}$ of pairs of lines ($1$-dimensional $K$-subspaces)
in $V$ such that
\begin{enumerate}
\item $\lambda_i^1 + \lambda_i^2$ is a hyperbolic plane for all $1 \leq i \leq
n$,
\item $\lambda_i^1 + \lambda_i^2$ is orthogonal to $\lambda_j^1 + \lambda_j^2$
for all $i \neq j$, and
\item $V = (\lambda_1^1 + \lambda_1^2) + \cdots + (\lambda_n^1 + \lambda_n^2)$.
\end{enumerate}
A vertex $t \in \Delta_n$ lies in the apartment specified by the frame
$\{\lambda_1^1, \lambda_1^2\}, \ldots,$ $\{\lambda_n^1, \lambda_n^2\}$ if for
any representative $L \in t$, there are lattices $M_i^j$ in $\lambda_i^j$ for
all $i, j$ such that $L = M_1^1 + M_1^2 + \cdots + M_n^1 + M_n^2$. The
following lemma is easily established.
\begin{nlem} \label{symplecticframe}
$\mbox{}$
\begin{enumerate}
\item Every symplectic basis for $V$ specifies an apartment of $\Delta_n$.
\item If $\Sigma$ is an apartment of $\Delta_n$, there is a symplectic basis
$\{u_1, \ldots, u_n,$ $w_1, \ldots, w_n\}$ for $V$ such that every vertex in
$\Sigma$ has the form
\[
[\cO\pi^{a_1}u_1 + \cdots + \cO\pi^{a_n}u_n + \cO\pi^{b_1}w_1 + \cdots +
\cO\pi^{b_n}w_n]
\]
for some $a_i, b_i \in \Z$.
\end{enumerate}
\end{nlem}

Since $\pi$ is fixed, if $\cB = \{u_1, \ldots, u_n, w_1, \ldots, w_n\}$ is a
symplectic basis for $V$, follow \cite[p.\ 3411]{shemanske} and write $(a_1,
\ldots, a_n; b_1, \ldots, b_n)_{\cB}$ for the lattice $\cO\pi^{a_1}u_1 + \cdots
+ \cO\pi^{a_n}u_n + \cO\pi^{b_1}w_1 + \cdots + \cO\pi^{b_n}w_n$ and $[a_1,
\ldots, a_n; b_1, \ldots, b_n]_{\cB}$ for its homothety class. Then the lattice
$L = (a_1, \ldots, a_n; b_1, \ldots, b_n)_{\cB}$ is primitive if and only if
$a_i + b_i = 0$ for all $i$ by \cite[p.\ 3411]{shemanske}, and $[L]$ is a
\emph{special} vertex in $\Delta_n$ if and only if $a_i + b_i = \mu$ is constant
for all $i$ by \cite[Corollary 3.4]{shemanske}. Note that by \cite[p.\
3412]{shemanske}, a chamber in $\Delta_n$ has exactly two special vertices.
\begin{nprop} \label{chambercount}
Every special vertex in $\Delta_n$ is contained in exactly
\[
\prod_{m = 1}^n\ \frac{q^{2m} - 1}{q - 1}
\]
chambers in $\Delta_n$.
\begin{proof}
Let $t \in \Delta_n$ be a special vertex. Then the number of chambers in
$\Delta_n$ containing $t$ is the number of chambers in the spherical $C_n(k)$
building (cf.\ \cite[p.\ 138]{brown}). By \cite[p.\ 6]{ronan}, a chamber in the
spherical $C_n(k)$ building is a maximal flag of non-trivial, \ti\ subspaces of
a $2n$-dimensional $k$-vector space endowed with a non-degenerate, alternating
bilinear form. An obvious modification of the proof of \cite[Proposition
2.4]{ss} finishes the proof.
\end{proof}
\end{nprop}

Let $C \in \Delta_n$ be a chamber with vertices $t_0, \ldots, t_n$, and let
$L_i \in t_i$ be representatives such that $L_0$ is primitive, $\langle L_i,
L_i \rangle \subseteq \pi\cO$ for all $1 \leq i \leq n$, and $\pi L_0
\subsetneq L_1 \subsetneq \cdots \subsetneq L_n \subsetneq L_0$. Let $\Sigma$
be an apartment of $\Delta_n$ containing $C$ and $\cB$ a symplectic basis for
$V$ specifying $\Sigma$ as in Lemma \ref{symplecticframe}. For all $0 \leq j
\leq n$, let
\[
L_j = (a_1^{(j)}, \ldots, a_n^{(j)}; b_1^{(j)}, \ldots, b_n^{(j)})_{\cB}.
\]
\begin{nlem} \label{chamber_special}
The two special vertices in $C$ are $[L_0]$ and $[L_n]$.
\begin{proof}
The fact that $[L_0]$ is special follows from \cite[Corollary 3.4]{shemanske}
and \cite[p.\ 3411]{shemanske}. To see that $[L_n]$ is special, note that if
$L_j$ represents a special vertex in $C$ for $1 \leq j \leq n$, then $a_i^{(j)}
+ b_i^{(j)} = \mu$ for all $i$ (by \cite[Corollary 3.4]{shemanske}), where $\mu
\in \{1, 2\}$ (since $\langle L_i, L_i \rangle \subseteq \pi\cO$). But $\mu =
2$ implies $L_j = \pi L_0$, which is impossible. Thus, $a_i^{(j)} + b_i^{(j)} =
1$ for all $i$ and $L_j/\pi L_0 \cong k^n$; hence, $j = n$.
\end{proof}
\end{nlem}

Let
\begin{align*}
\GSp_n(K) &= \{g = \left(\begin{smallmatrix}
A & B \\
C & D
\end{smallmatrix}\right) \in M_{2n}(K) : A, B, C, D \in M_n(K),\ A^tC = C^tA,\\
&\qquad B^tD = D^tB,\ A^tD - C^tB = \nu(g)I_n \text{ for some $\nu(g) \in
K^\times$}\};
\end{align*}
alternatively, abuse notation and think of $\GSp_n(K)$ as
\begin{align*}
\{g \in \GL_K(V) &: \text{$\forall\ v_1, v_2 \in V$, $\exists\ \nu(g) \in
K^\times$ such that} \\
& \qquad \langle gv_1, gv_2 \rangle = \nu(g)\langle v_1, v_2 \rangle\}.
\end{align*}
Then $\Sp_n(K)$ consists of those $g \in \GSp_n(K)$ such that $\nu(g) = 1$;
alternatively, $\langle gv_1, gv_2 \rangle = \langle v_1, v_2 \rangle$ for all
$v_1, v_2 \in V$. Let $\cB = \{u_1, \ldots, u_n, w_1, \ldots,$ $w_n\}$ be a
symplectic basis for $V$ and $g \in \GSp_n(K)$. If $t \in \Delta_n$ is a vertex
with representative $L = (a_1, \ldots, a_n; b_1, \ldots, b_n)_{\cB}$, define
\[
gt = [\cO\pi^{a_1}gu_1 + \cdots + \cO\pi^{a_n}gu_n + \cO\pi^{b_1}gw_1 + \cdots
+ \cO\pi^{b_n}gw_n].
\]
Note that
\[
\cB_g := \{\nu(g)^{-1}gu_1, \ldots, \nu(g)^{-1}gu_n, gw_1, \ldots, gw_n\}
\]
is a symplectic basis for $V$; hence, $m = \ord(\nu(g))$ implies $gL = (a_1 +
m, \ldots, a_n + m; b_1, \ldots, b_n)_{\cB_g}$.
\begin{nprop} \label{GSpaction}
The group $\GSp_n(K)$ acts transitively on the special vertices in $\Delta_n$.
\begin{proof}
Note that if $\GSp_n(K)$ acts on the special vertices in $\Delta_n$, then
\cite[Proposition 3.3]{shemanske} implies that the action is transitive. We
thus show that $\GSp_n(K)$ acts on the special vertices in $\Delta_n$. Let $t
\in \Delta_n$ be a special vertex and $L \in t$ a representative such that
there is a primitive lattice $L_0$ with $\langle L, L \rangle \subseteq \pi\cO$
and $\pi L_0 \subseteq L \subseteq L_0$. Let $\Sigma$ be an apartment of
$\Delta_n$ containing $t$ and $[L_0]$, and let $\cB$ be a symplectic basis for
$V$ specifying $\Sigma$ as in Lemma \ref{symplecticframe}. Then \cite[p.\
3411]{shemanske}, the last lemma, and \cite[Corollary 3.4]{shemanske} imply
\[
L_0 = (c_1, \ldots, c_n; -c_1, \ldots, -c_n)_{\cB}
\]
and
\[
L = (a_1, \ldots, a_n; \mu - a_1, \ldots, \mu - a_n)_{\cB},
\]
where $\mu \in \{1, 2\}$. Let $g \in \GSp_n(K)$ with $\ord(\nu(g)) = m$. Since
$gt = [a_1 + m, \ldots, a_n + m; \mu - a_1, \ldots, \mu - a_n]_{\cB_g}$,
\cite[Corollary 3.4]{shemanske} implies that it suffices to show $gt$ is a
vertex in $\Delta_n$. First suppose $m = 2r$ for some $r \in \Z$. Then
$\pi^{-r}gL_0$ is primitive, $\langle \pi^{-r}gL, \pi^{-r}gL \rangle \subseteq
\pi\cO$, and $\pi^{-r}g(\pi L_0) \subseteq \pi^{-r}gL \subseteq \pi^{-r}gL_0$;
i.e., $gt$ is a vertex in $\Delta_n$. Now suppose $m = 2r + 1$. If $\mu = 1$,
then $\pi^{-r - 1}gL$ is primitive and $gt$ is a vertex in $\Delta_n$.
Otherwise, $\mu = 2$, and $\langle \pi^{-r - 1}gL, \pi^{-r - 1}gL \rangle
\subseteq \pi\cO$. Let $\pi M_0 = (a_1 + r, \ldots, a_n + r; \mu - a_1 - r,
\ldots, \mu - a_n - r)_{\cB_g}$. Then $M_0$ is primitive and $\pi M_0 \subseteq
\pi^{-r - 1}gL \subseteq M_0$; i.e., $gt$ is a vertex in $\Delta_n$. Thus,
$\GSp_n(K)$ acts on the special vertices in $\Delta_n$.
\end{proof}
\end{nprop}

Call two distinct, incident vertices in $\Delta_n$ \emph{adjacent}.
\begin{nprop} \label{adjacent_special}
The group $\GSp_n(K)$ takes adjacent special vertices in $\Delta_n$ to adjacent
special vertices in $\Delta_n$.
\begin{proof}
Let $t, t' \in \Delta_n$ be adjacent special vertices, and let $L \in t$ and
$L' \in t'$ be representatives such that there is a primitive lattice $L_0$
with $\langle L, L \rangle \subseteq \pi\cO$, $\langle L', L' \rangle \subseteq
\pi\cO$, and either $\pi L_0 \subseteq L \subsetneq L' \subseteq L_0$ or $\pi
L_0 \subseteq L' \subsetneq L \subseteq L_0$. By Lemma \ref{chamber_special},
either $L' = L_0$ or $L = L_0$; i.e., either $\pi L' \subsetneq L \subsetneq
L'$ with $L'$ primitive or $\pi L \subsetneq L' \subsetneq L$ with $L$
primitive. By symmetry, assume $\pi L \subsetneq L' \subsetneq L$ with $L$
primitive, and let $g \in \GSp_n(K)$ with $\ord(\nu(g)) = m$. If $m = 2r$ for
some $r \in \Z$, then (as in the last proposition) $\pi^{-r}gL$ is primitive,
$\langle \pi^{-r}gL', \pi^{-r}gL' \rangle \subseteq \pi\cO$, and $\pi^{-r}g(\pi
L) \subsetneq \pi^{-r}gL' \subsetneq \pi^{-r}gL$. Similarly, if $m = 2r + 1$,
then $\pi^{-r - 1}gL'$ is primitive, $\langle \pi^{-r}gL, \pi^{-r}gL \rangle
\subseteq \pi\cO$, and $\pi(\pi^{-r - 1}gL') \subsetneq \pi^{-r}gL \subsetneq
\pi^{-r - 1}gL'$.
\end{proof}
\end{nprop}

\begin{nlem}
If $t, t' \in \Delta_n$ are adjacent special vertices and $t$ has a primitive
representative $L$, then there is a representative $L' \in t'$ with $\pi L
\subsetneq L' \subsetneq L$ such that the number of chambers in $\Delta_n$
containing both $t$ and $t'$ equals the number of maximal flags of non-trivial,
proper $k$-subspaces of $L'/\pi L$.
\begin{proof}
Since a chamber $C \in \Delta_n$ containing both $t$ and $t'$ has $n + 1$
vertices $t_0, \ldots, t_n$ that have representatives $L_i \in t_i$ such that
$L_0$ is primitive, $\langle L_i, L_i \rangle \subseteq \pi\cO$ for all $1 \leq
i \leq n$, and $\pi L_0 \subsetneq L_1 \subsetneq \cdots \subsetneq L_n
\subsetneq L_0$, Lemma \ref{chamber_special} implies $L_0$ and $L_n$ represent
the two special vertices in $C$; in particular, $t = t_0$ and $t' = t_n$. Let
$L = L_0$ and $L' = L_n$. Then $\pi L \subsetneq L' \subsetneq L$ and varying
$L_1, \ldots, L_{n - 1}$ over all lattices in $V$ contained in $L$ and
containing $\pi L$ such that $\langle L_i, L_i \rangle \subseteq \pi\cO$ for
all $1 \leq i \leq n - 1$ and $\pi L \subsetneq L_1 \subsetneq \cdots
\subsetneq L_{n - 1} \subsetneq L'$ gives all the chambers in $\Delta_n$
containing both $t$ and $t'$.
\end{proof}
\end{nlem}
\begin{nprop} \label{adjacent}
If $t \in \Delta_n$ is a special vertex, then $t$ is adjacent to exactly
$\prod_{m = 1}^n\ (q^m + 1)$ distinct special vertices in $\Delta_n$.
\begin{proof}
Let $t \in \Delta_n$ be a special vertex. By Proposition \ref{chambercount},
the number of chambers in $\Delta_n$ containing $t$ is $\prod_{m = 1}^n\
((q^{2m} - 1)/(q - 1))$. Since this counts a special vertex $t' \in \Delta_n$
adjacent to $t$ more than once if there is more than one chamber in $\Delta_n$
containing both $t$ and $t'$, the last lemma and \cite[Proposition 2.4]{ss}
finish the proof.
\end{proof}
\end{nprop}
\begin{nprop} \label{GSpsymplecticstab}
Let $L$ be a lattice in $V$. Then $\GSp_n(\cO) := \GSp_n(K) \cap \GL_{2n}(\cO)$
can be identified with $\{g \in \GSp_n(K) : gL = L\}$, where $g$ acts on $L$ as
the matrix of a linear transformation with respect to a fixed basis for $L$.
\begin{proof}
It suffices to prove the proposition for any lattice $L$ in $V$. Let $\cB =
\{e_1, \ldots, e_{2n}\}$ be the standard unit basis for $V$, and let $L = (0,
\ldots, 0; 0, \ldots,$ $0)_{\cB}$. Let $g \in \GSp_n(K)$ such that $gL = L$. Then
$g \in \GL_{2n}(\cO)$; i.e., $g \in \GSp_n(K) \cap \GL_{2n}(\cO)$. On the other
hand, any element of $\GSp_n(K) \cap \GL_{2n}(\cO)$ fixes $L$.
\end{proof}
\end{nprop}

We end this section with some terminology used in Proposition \ref{connected}.
By \cite[p.\ 13]{brown}, two chambers in $\Delta_n$ are \emph{adjacent} if they
have a common codimension-one face, and a \emph{gallery} in $\Delta_n$
connecting the chambers $C, C' \in \Delta_n$ is a sequence $C = C_0, \ldots,
C_m = C'$ of chambers in $\Delta_n$ such that $C_i$ and $C_{i + 1}$ are
adjacent for all $0 \leq i \leq m - 1$. Recall that any two chambers in a
building can be connected by a gallery.
\section{The Isoperimetric Constant of $Y_n$} \label{comp}

The graph theory notation and terminology used in this section is primarily
from \cite[p.\ 7]{woess}; any terminology not defined there is from \cite[pp.\
1 -- 4]{bollobas}, except for the definition of walk, which is from \cite[p.\
2]{murty}. In particular, a graph is a finite or countably infinite set $X$ of
vertices, together with a symmetric neighborhood or adjacency relation $\sim$.
Let $X$ be a connected, $r$-regular ($r$ finite) graph with infinitely many
vertices. As in \cite[p.\ 116]{bms}, if $X'$ is a subset of $X$, let $\partial
X'$ denote the set of edges in $X$ incident to exactly one vertex in $X'$. Then
the \emph{isoperimetric constant} of $X$ is $h(X) :=
\inf(\Card(\partial X')/\Card(X'))$, where the infimum is over all finite,
non-empty subsets $X'$ of $X$. Note that $h(X)$ is related to the spectral
radius of $X$: recall that the adjacency operator $A(X)$ of $X$ acts on the
Hilbert space of functions $f: X \rightarrow \C$ such that $\sum_{x \in X}\
|f(x)|^2 < \infty$. Then the spectrum of $A(X)$ is $\{\lambda \in \C : A(X) -
\lambda I\text{ is not invertible}\}$, and the \emph{spectral radius} of $X$ is
\[
\rho(X) := \sup\{|\lambda| : \lambda \text{ is in the spectrum of $A(X)$}\};
\]
equivalently, $\rho(X) = \|A(X)\|$, the norm of $A(X)$, by \cite[p.\
252]{mohar} ($A(X)$ is bounded by \cite[Theorem 3.2]{mohar} and can be shown to
be self-adjoint). Then by \cite[Theorem 3.1]{bms}, for a connected, $r$-regular
graph $X$ with infinitely many vertices, to show $h(X) > 0$, it suffices to
compute the spectral radius $\rho(X)$ of $X$ and show $\rho(X) < r$.

Let $Y_n$ be the subgraph of $\Delta_n$ induced by the \emph{special} vertices
in $\Delta_n$. Then by Proposition \ref{adjacent}, $Y_n$ is $(\prod_{m = 1}^n\
(q^m + 1))$-regular. Moreover, \cite[p.\ 3414]{shemanske} implies that $Y_n$
has infinitely many vertices.
\begin{nprop} \label{connected}
The graph $Y_n$ is connected.
\begin{proof}
Let $t \neq t' \in Y_n$. Since there is nothing to prove if $t$ and $t'$ are in
a common chamber in $\Delta_n$, assume that no chamber in $\Delta_n$ contains
both $t$ and $t'$. Let $C, C' \in \Delta_n$ be chambers such that $t \in C$ and
$t' \in C'$, and let $C = C_0, \ldots, C_m = C'$ be a gallery in $\Delta_n$
connecting $C$ and $C'$. For all $0 \leq i \leq m - 1$, let $S_i$ be the set of
special vertices contained in $C_i$ and $C_{i + 1}$. Let $t_0 = t$, $t_{m + 1}
= t'$, and $t_i \in S_{i - 1}$ for all $1 \leq i \leq m$. Then $t_i$ and $t_{i
+ 1}$ are incident vertices in $\Delta_n$ for all $i$; hence, for each $0 \leq
i \leq m$, either $t_i = t_{i + 1}$ or $t_i$ and $t_{i + 1}$ are adjacent
vertices in $Y_n$. But we can remove vertices from the sequence $t_0, \ldots,
t_{m + 1}$ until we are left with a walk in $Y_n$ connecting $t$ and $t'$
(since $Y_n$ has no loops, we need successive vertices to be adjacent).
\end{proof}
\end{nprop}

If $G$ is a group acting on a set $S$ and $a, b \in S$, write $G_a$ for
$\Stab_G(a) = \{g \in G : ga = a\}$ and $G_ab$ for $\{gb : g \in G_a\}$. The
main tool that we use is the following reformulation of a special case of
\cite[Theorem (12.10)]{woess}.
\begin{nthm} \label{thm}
Let $X$ be a locally finite, regular graph. If there is a solvable group $Q$
that acts transitively on $X$ as a group of automorphisms and has a left Haar
measure $\mu(\cdot)$, then the spectral radius $\rho(X)$ of $X$ is
\begin{equation} \label{formula}
\rho(X) = \sum_{t' \sim t_0}\ \sqrt{\frac{\Card(Q_{t'}t_0)}{\Card(Q_{t_0}t')}},
\end{equation}
where $t_0$ is any vertex in $X$.
\end{nthm}

Recall that $\PGSp_n(K) = \GSp_n(K)/K^\times$, where we identify $K^\times$
with the scalar matrices of $\GSp_n(K)$. As in \cite[Example (12.20)]{woess},
write $g$ as a matrix in $\GSp_n(K)$ while thinking of it as an element of
$\PGSp_n(K)$ consisting of all its non-zero multiples. Let $Q$ be the image in
$\PGSp_n(K)$ of $Q'$, where
\[
Q' = \left\{
\left(\begin{smallmatrix}
A & B \\
0 & D
\end{smallmatrix}\right) \in \GSp_n(K) : \text{$A \in M_n(K)$ is upper
triangular}\right\},
\]
a minimal parabolic subgroup of $\GSp_n(K)$, and let $V$ be a $2n$-dimensional
$K$-vector space with a non-degenerate, alternating bilinear form $\langle
\cdot, \cdot \rangle$. We want to apply Theorem \ref{thm} to the graph $Y_n$
with $Q$ as above and $t_0 = o = [\cO^{2n}] = [0, \ldots, 0; 0, \ldots,
0]_{\cB_0}$, where $\cB_0:= \{e_1, \ldots, e_n, f_1, \ldots, f_n\}$ is the
standard symplectic basis for $V$ ($f_i = e_{n + i}$ for all $i$). Note that a
modification of the proof of Proposition 1.3.7 of \cite{andrianov} shows that
for any $h \in \GSp_n(K)$, there is a $g \in \Sp_n(\cO)$ such that $gh^{-1} \in
Q'$; furthermore, $(hg^{-1})o = ho$ by Proposition \ref{GSpsymplecticstab}.
This, together with Proposition \ref{GSpaction}, implies that $Q$ acts
transitively on the vertices in $Y_n$. The fact that $Q$ acts on $Y_n$ as a
group of automorphisms follows from the fact that $Q'$ does (see Proposition
\ref{adjacent_special}). Since the group of upper triangular matrices in
$\GL_{2n}(K)$ is solvable, so is $Q$. Verifying that $Q$ has a left Haar
measure is a straightforward exercise involving topological groups.

\begin{figure}
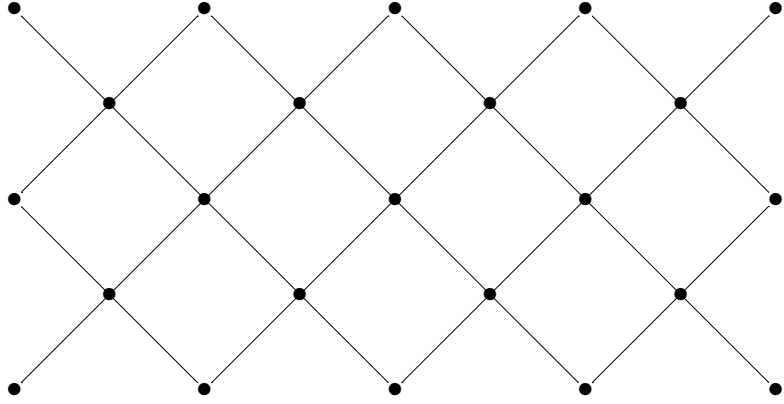

\[
\xy 0;/r0.20pc/:
(-60, -30)*{\bullet}="dd0";
(-30, -30)*{\bullet}="dd1";
(0, -30)*{\bullet}="dd2";
(30, -30)*{\bullet}="dd3";
(60, -30)*{\bullet}="dd4";
(-45, -15)*{\bullet}="d0";
(-15, -15)*{\bullet}="d1";
(15, -15)*{\bullet}="d2";
(45, -15)*{\bullet}="d3";
(-60, 0)*{\bullet}="c0";
(-30, 0)*{\bullet}="c1";
(0, 0)*{\bullet}="c2";
(30, 0)*{\bullet}="c3";
(60, 0)*{\bullet}="c4";
(-45, 15)*{\bullet}="u0";
(-15, 15)*{\bullet}="u1";
(15, 15)*{\bullet}="u2";
(45, 15)*{\bullet}="u3";
(-60, 30)*{\bullet}="uu0";
(-30, 30)*{\bullet}="uu1";
(0, 30)*{\bullet}="uu2";
(30, 30)*{\bullet}="uu3";
(60, 30)*{\bullet}="uu4";
{\ar@{-} "dd0";"uu2"};
{\ar@{-} "dd1";"c0"};
{\ar@{-} "dd1";"uu3"};
{\ar@{-} "dd2";"uu0"};
{\ar@{-} "dd2";"uu4"};
{\ar@{-} "dd3";"uu1"};
{\ar@{-} "dd3";"c4"};
{\ar@{-} "dd4";"uu2"};
{\ar@{-} "c0";"uu1"};
{\ar@{-} "c4";"uu3"};
\endxy
\]
\caption{Subgraph of $Y_2$.}
\label{Sp2apt_special}
\end{figure}

Since we take $t_0 = o$ in \eqref{formula}, it suffices to determine
$\Card(Q_{t'}o)$ and $\Card(Q_ot')$ for vertices $t' \in Y_n$ with $t' \sim o$.
For a symplectic basis $\cB$ for $V$ and $\Sigma$ the apartment of $\Delta_n$
specified by $\cB$ as in Lemma \ref{symplecticframe}, follow \cite[Example
(12.20)]{woess} and write $U(\cB)$ for the subgraph of $Y_n$ induced by the
special vertices in $\Sigma$ (see Figure \ref{Sp2apt_special} for a partial
picture of $U(\cB)$ when $n = 2$).
Let $\cE(n) = \{0, 1\}^n$. Then the neighbors of $o$ in $U(\cB_0)$ are
\[
x_{\ue} := [\varepsilon_1, \ldots, \varepsilon_n; 1 - \varepsilon_1, \ldots, 1
- \varepsilon_n]_{\cB_0},
\]
where $\ue = (\varepsilon_1, \ldots, \varepsilon_n) \in \cE(n)$, and
Proposition \ref{GSpsymplecticstab} implies that the stabilizer of $o$ in
$\PGSp_n(K)$ can be identified with $\PGSp_n(\cO)$. It follows that if $g =
(g_{ij}) \in Q_o = Q \cap \PGSp_n(\cO)$ and $|\cdot|_K$ is the absolute value
of $K$, normalized such that $|\pi|_K = 1/q$, then
\begin{equation} \label{absstaboentries}
|g_{ii}|_K = 1, \qquad |g_{ij}|_K \leq 1 \text{ if $i < j$ (with $1 \leq i \leq
n$) or $2n \geq i > j \geq n + 1$},
\end{equation}
and $g_{ij} = 0$ otherwise. For $\ue \in \cE(n)$, let
\[
g_{\ue} := \diag(\pi^{\varepsilon_1}, \ldots, \pi^{\varepsilon_n}, \pi^{1 -
\varepsilon_1}, \ldots, \pi^{1 - \varepsilon_n}),
\]
and note that $g_{\ue}o = x_{\ue}$ and $Q_{x_{\ue}} = g_{\ue}Q_og_{\ue}^{-1}$.
Then \eqref{absstaboentries} and a similar analysis for an element of
$Q_{x_{\ue}}$ imply that for $h = (h_{ij}) \in Q_o \cap Q_{x_{\ue}}$,
\begin{align}
|h_{ij}|_K &= 1 \text{ if $i = j$}, \label{intersectdiagonal} \\
|h_{ij}|_K &\leq \begin{cases}
q^{-\max\{0, \varepsilon_i - \varepsilon_j\}} & \text{if $1 \leq i < j \leq
n$}, \\
q^{-\max\{0, \varepsilon_i - (1 - \varepsilon_{j - n})\}} & \text{if $1 \leq i
\leq n < j \leq 2n$}, \\
q^{-\max\{0, -\varepsilon_{i - n} + \varepsilon_{j - n}\}} & \text{if $n + 1
\leq j < i \leq 2n$},
\end{cases} \label{intersectother}
\end{align}
and $h_{ij} = 0$ otherwise. In addition, by the orbit-stabilizer theorem,
\[
\Card(Q_ox_{\ue}) = \frac{\mu(Q_o)}{\mu(Q_o \cap Q_{x_{\ue}})}
\qquad \text{and} \qquad \Card(Q_{x_{\ue}}o) = \Card(Q_ox_{\underline{1} -
\ue}),
\]
where $\mu(\cdot)$ is a left Haar measure on $Q$ and $\underline{1} = (1,
\ldots, 1) \in \cE(n)$.

Consequently, to determine $\Card(Q_ox_{\ue})$ and $\Card(Q_{x_{\ue}}o)$, it
suffices to find $\mu(Q_o)$ and $\mu(Q_o \cap Q_{x_{\ue}})$. Since we can
completely characterize both $Q_o$ and $Q_o \cap Q_{x_{\ue}}$ in terms of
$|\cdot|_K$, $\mu(Q_o)$ (resp., $\mu(Q_o \cap Q_{x_{\ue}})$) is the product
of the Haar measures of the unconstrained, non-constant entries of an element
$g \in Q_o$ (resp., of $g \in Q_o \cap Q_{x_{\ue}}$). Write $\vol(g_{ij})$ for
the Haar measure of the $(i, j)$-entry of $g = (g_{ij})$, with $\vol$
normalized such that $\vol(\cO) = 1$, $\vol(\pi\cO) = 1/q$, and $\vol(\cO^\times)
= 1$. It follows from \eqref{absstaboentries} that
\[
\mu(Q_o) = 1.
\]

For $\ue \in \cE(n)$, let $|\ue|_n := \sum_{i = 1}^n \varepsilon_i$. Then
\[
M(\ue, n) := \sum_{1 \leq i < j \leq n}\ \max\{0, \varepsilon_i -
\varepsilon_j\}
\]
counts the number of $0$s that follow each $1$ in $\ue$. Moreover,
\[
\sum_{1 \leq i \leq n < j \leq i + n}\ \max\{0, \varepsilon_i - (1 -
\varepsilon_{j - n})\} = \sum_{1 \leq j \leq i \leq n}\ \max\{0, (\varepsilon_i
+ \varepsilon_j) - 1\}
\]
adds the number of $1$s in $\ue$ to the number of $1$s that follow each $1$ in
$\ue$; i.e., if $|\ue|_n = m$, then the last sum is $m(m + 1)/2$. Let $\ue \in
\cE(n)$ with $|\ue|_n = m$ and $h = (h_{ij}) \in Q_o \cap Q_{x_{\ue}}$. Then
\eqref{intersectdiagonal} and \eqref{intersectother} imply
\[
\mu(Q_o \cap Q_{x_{\ue}}) = q^{-M(\ue, n)}q^{-m(m + 1)/2};
\]
hence,
\begin{equation} \label{orbitmeasure}
\Card(Q_ox_{\ue}) = q^{M(\ue, n)}q^{m(m + 1)/2}
\end{equation}
and
\[
\Card(Q_{x_{\ue}}o) = q^{M(\underline{1} - \ue, n)}q^{(n - m)(n - m + 1)/2}.
\]
Since $M(\underline{1} - \ue, n)$ counts the number of $0$s that precede each
$1$ in $\ue$, $M(\ue, n) + M(\underline{1} - \ue, n)$ is the product of the
number of zeros in $\ue$ and $|\ue|_n$. Consequently,
\begin{equation} \label{product}
\Card(Q_ox_{\ue})\Card(Q_{x_{\ue}}o) = q^{n(n + 1)/2}.
\end{equation}

\begin{nlem}
Let $\ue \neq \ue' \in \cE(n)$. Then for all $g \in Q_o$, $x_{\ue'} \neq
gx_{\ue}$.
\begin{proof}
Let $\ue = (\varepsilon_1, \ldots, \varepsilon_n)$ and $\ue' = (\varepsilon_1',
\ldots, \varepsilon_n')$, and let $1 \leq \ell \leq n$ such that
$\varepsilon_\ell \neq \varepsilon_\ell'$. If $x_{\ue'} = gx_{\ue}$ for some $g
= (g_{ij}) \in Q_o$, then $g_{\ue}g_{\ue'}^{-1}g \in Q_{x_{\ue}}$, which is
impossible since by \eqref{absstaboentries}, the $(\ell, \ell)$-entry of
$g_{\ue}g_{\ue'}^{-1}g$ satisfies
\[
\Bigl|g_{\ell\ell}\pi^{\varepsilon_\ell - \varepsilon_\ell'}\Bigr|_K \in
\{q^{-1}, q\},
\]
contradicting \eqref{intersectdiagonal}.
\end{proof}
\end{nlem}

It follows that for $\ue, \ue' \in \cE(n)$, the sets $Q_ox_{\ue}$ and
$Q_ox_{\ue'}$ are disjoint if $\ue \neq \ue'$.
\begin{nprop} \label{oneighbors}
The set of vertices in $Y_n$ adjacent to $o$ is $N(o) = \cup_{\ue \in \cE(n)}\
Q_ox_{\ue}$.
\begin{proof}
First note that by the last lemma and \eqref{orbitmeasure},
\[
\Card\left(\bigcup_{\ue \in \cE(n)}\ Q_ox_{\ue}\right) = \sum_{m = 0}^n\ q^{m(m
+ 1)/2}\left(\sum_{\substack{\ue \in \cE(n),\\
|\ue|_n = m}}\ q^{M(\ue, n)}\right).
\]
Let
\[
W(n, m) = \sum_{\substack{\ue \in \cE(n),\\ |\ue|_n = m}}\ q^{M(\ue, n)},
\]
and note that by \cite[p.\ 135]{woess}, $W(n, m) = {n \choose m}_q$, the number
of $m$-dimensional subspaces of an $n$-dimensional $\F_q$-vector space, for all
$0 \leq m \leq n$ (see \cite[p.\ 133]{woess} for the formula for ${n \choose
m}_q$). The proof now follows from Proposition \ref{adjacent} since for all $n
\geq 1$,
\[
\prod_{m = 1}^n\ (1 + q^m) = \sum_{i = 0}^n\ q^{i(i + 1)/2}{n \choose i}_q.
\qedhere
\]
\end{proof}
\end{nprop}
\begin{ncor}
If $t \in \Delta_n$ is a vertex with a primitive representative, then the set
of vertices in $Y_n$ adjacent to $t$ is $\cup_{\ue \in \cE(n)}\ Q_tx_{\ue}$.
\begin{proof}
This follows from the last proposition since we only used the fact that $o$ has
a primitive representative.
\end{proof}
\end{ncor}

\begin{nprop} \label{Spnspecrad}
The spectral radius of $Y_n$ is $\rho(Y_n) = 2^nq^{n(n + 1)/4}$.
\begin{proof}
This follows from \eqref{formula}, Proposition \ref{oneighbors},
\eqref{product}, and the fact that
\[
\frac{\Card(Q_{x_{\ue}}o)}{\Card(Q_ox_{\ue})} =
\frac{\Card(Q_{gx_{\ue}}o)}{\Card(Q_ogx_{\ue})}
\]
for all $g \in Q_o$.
\end{proof}
\end{nprop}
\begin{nthm} \label{specrad}
The isoperimetric constant of $Y_n$ satisfies $h(Y_n) > 0$.
\begin{proof}
By \cite[Theorem 3.1]{bms} and the last proposition, it suffices to show that
$2^nq^{n(n + 1)/4} < \prod_{m = 1}^n\ (q^m + 1)$, which is straightforward.
\end{proof}
\end{nthm}

\begin{figure}
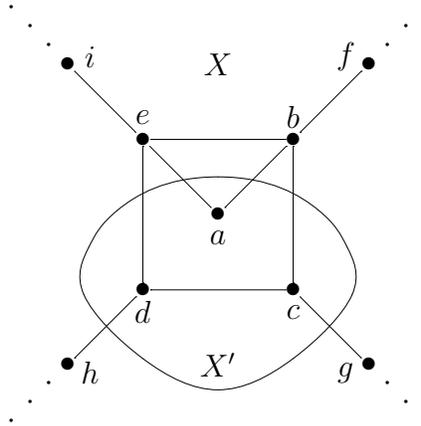

\centering
\[
{\xy
(20, -20)*{\bullet}="2";
(17, -21)*{\text{$g$}};
(22.5, -22.5)*{\cdot};
(25, -25)*{\cdot};
(27.5, -27.5)*{\cdot};
(20, 20)*{\bullet}="1";
(17, 21)*{\text{$f$}};
(22.5, 22.5)*{\cdot};
(25, 25)*{\cdot};
(27.5, 27.5)*{\cdot};
(10, -10)*{\bullet}="c";
(10, -13)*{\text{$c$}};
(10, 10)*{\bullet}="b";
(10, 13)*{\text{$b$}};
(0, 0)*{\bullet}="a";
(0, -3)*{\text{$a$}};
(-10, -10)*{\bullet}="d";
(-10, -13)*{\text{$d$}};
(-10, 10)*{\bullet}="e";
(-10, 13)*{\text{$e$}};
(-20, -20)*{\bullet}="3";
(-17, -21)*{\text{$h$}};
(-22.5, -22.5)*{\cdot};
(-25, -25)*{\cdot};
(-27.5, -27.5)*{\cdot};
(-20, 20)*{\bullet}="4";
(-17, 21)*{\text{$i$}};
(-22.5, 22.5)*{\cdot};
(-25, 25)*{\cdot};
(-27.5, 27.5)*{\cdot};
{\ar@{-} "a";"e" };
{\ar@{-} "a";"b" };
{\ar@{-} "e";"b" };
{\ar@{-} "e";"d" };
{\ar@{-} "b";"c" };
{\ar@{-} "d";"c" };
{\ar@{-} "b";"1" };
{\ar@{-} "e";"4" };
{\ar@{-} "d";"3" };
{\ar@{-} "c";"2" };
(0, 5)*{}="A";
"A";"A"**\crv{(-7, 5) & (-15, 0) & (-17, -4) & (-20, -10) & (-10, -20)
& (0, -25) & (10, -20) & (20, -10) & (17, -4) & (15, 0) & (7, 5)};
(0, 20)*{\text{$X$}};
(0, -20)*{\text{$X'$}};
\endxy}
\]
\caption{In general, $\Card(\sigma(X')) \leq \Card(\partial X')$.}
\label{boundaryinequality}
\end{figure}

Let $X$ be a connected graph with infinitely many vertices and of bounded
degree. Following \cite[p.\ 2480]{et}, for a connected, induced subgraph $X'$
of $X$ with at least one edge, let $\sigma(X')$ denote the set of vertices in
$X'$ adjacent to a vertex in $X$ not in $X'$. Then $X$ is \emph{amenable} if
$\inf(\Card(\sigma(X'))/\Card(X'))$ $= 0$, where the infimum is over all finite,
connected, induced subgraphs $X'$ of $X$ with at least one edge. Note that if
$X$ is finite and has at least one edge, then $X$ is trivally amenable.
Furthermore, if $X'$ is a finite, connected, induced subgraph of $X$ with at
least one edge, then $\Card(\sigma(X')) \leq \Card(\partial X')$ (see Figure
\ref{boundaryinequality}).
On the other hand, since $X$ has bounded degree, $X$ amenable implies $h(X) =
0$.
\begin{ncor}
The graph $Y_n$ is non-amenable.
\begin{proof}
As noted above, $X$ amenable implies $h(X) = 0$; hence, $h(Y_n) > 0$ implies
$Y_n$ non-amenable.
\end{proof}
\end{ncor}

For $x \in Y_n$, let $B_i(x) = \{y \in Y_n : \dist(x, y) \leq i\}$, where
$\dist(x, y)$ is the (graph) distance between $x$ and $y$.
\begin{ncor} \label{constant}
There is a constant $C > 1$ such that for all $x \in Y_n$ and for all $i \in
\Z^{\geq 0}$, $\Card(B_i(x)) > C^i$.
\begin{proof}
This follows from the last theorem and \cite[Theorem 2.2]{bms}.
\end{proof}
\end{ncor}

Finally, note that the building $\Delta_n$ is a subcomplex of $\Xi_{2n}$
(compare the description of $\Delta_n$ given in Section \ref{prelims} with the
description of $\Xi_n$ in \cite[p.\ 115]{ronan}); hence, $Y_n$ is a subgraph of
the one-complex $X_{2n}$ of $\Xi_{2n}$. Since both $Y_n$ and $X_{2n}$ are
expanding, it is natural to ask about their relative expansion properties. It
is straightfoward to show that $\rho(Y_n) < \rho(X_{2n})$ for all $n \geq 2$,
but since the degree of any vertex in $Y_n$ is also strictly less than the
degree of any vertex in $X_{2n}$, it is unclear what this reveals. On the other
hand, the analogue of Corollary \ref{constant} holds for $X_{2n}$. Then
Theorems 2.2 and 3.1 of \cite{bms} provide a constant $C(Y_n)$ (resp.,
$C(X_{2n})$) as in Corollary \ref{constant} (resp., as in the analogue of
Corollary \ref{constant} for $X_{2n}$) in terms of $\rho(Y_n)$ and the degree
of any vertex in $Y_n$ (resp., in terms of $\rho(X_{2n})$ and the degree of any
vertex in $X_{2n}$). This thus raises the question of whether there is any
relationship between $C(Y_n)$ and $C(X_{2n})$. For example, our data for $2
\leq n \leq 5$ and $q = p^i$, where $1 \leq i \leq 5$ and $p$ is one of the
first five primes, indicates that $C(Y_n) < C(X_{2n})$; is this always the
case?
\bibliographystyle{abbrv}
\bibliography{blah}

\begin{thebibliography}{10}

\bibitem{andrianov}
A.~Andrianov.
\newblock {\em Quadratic Forms and {H}ecke Operators}, volume 286 of {\em
  Grundlehren der Mathematischen Wissenschaften}.
\newblock Springer-Verlag, Berlin, 1987.

\bibitem{bien}
F.~Bien.
\newblock Constructions of telephone networks by group representations.
\newblock {\em Notices of the American Mathematical Society}, 36(1):5 -- 22,
  1989.

\bibitem{bms}
N.~Biggs, B.~Mohar, and J.~Shawe-Taylor.
\newblock The spectral radius of infinite graphs.
\newblock {\em The Bulletin of the London Mathematical Society}, 20(2):116 --
  120, 1988.

\bibitem{bollobas}
B.~Bollob\'{a}s.
\newblock {\em Graph Theory: An Introductory Course}, volume~63 of {\em
  Graduate Texts in Mathematics}.
\newblock Springer-Verlag, New York-Berlin, 1979.

\bibitem{brown}
K.~Brown.
\newblock {\em Buildings}.
\newblock Springer-Verlag, New York, 1989.

\bibitem{csz}
D.~Cartwright, P.~Sol\'{e}, and A.~\.{Z}uk.
\newblock Ramanujan geometries of type $\widetilde{A}_n$.
\newblock {\em Discrete Mathematics}, 269(1 -- 3):35 -- 43, 2003.

\bibitem{dsv}
G.~Davidoff, P.~Sarnak, and A.~Valette.
\newblock {\em Elementary Number Theory, Group Theory, and {R}amanujan Graphs},
  volume~55 of {\em London Mathematical Society Student Texts}.
\newblock Cambridge University Press, Cambridge, 2003.

\bibitem{et}
G.~Elek and G.~Tardos.
\newblock On roughly transitive amenable graphs and harmonic {D}irichlet
  functions.
\newblock {\em Proceedings of the American Mathematical Society}, 128(8):2479
  -- 2485, 2000.

\bibitem{garrett}
P.~Garrett.
\newblock {\em Buildings and Classical Groups}.
\newblock Chapman \& Hall, London, 1997.

\bibitem{mohar}
B.~Mohar.
\newblock The spectrum of an infinite graph.
\newblock {\em Linear Algebra and its Applications}, 48:245 -- 256, 1982.

\bibitem{murty}
M.~Murty.
\newblock Ramanujan graphs.
\newblock {\em Journal of the Ramanujan Mathematical Society}, 18(1):33 -- 52,
  2003.

\bibitem{parkinsonmz}
J.~Parkinson.
\newblock Spherical harmonic analysis on affine buildings.
\newblock {\em Mathematische Zeitschrift}, 253(3):571 -- 606, 2006.

\bibitem{ronan}
M.~Ronan.
\newblock {\em Lectures on Buildings}, volume~7 of {\em Perspectives in
  Mathematics}.
\newblock Academic Press, Inc., Boston, MA, 1989.

\bibitem{sw}
L.~Saloff-Coste and W.~Woess.
\newblock Transition operators, groups, norms, and spectral radii.
\newblock {\em Pacific Journal of Mathematics}, 180(2):333 -- 367, 1997.

\bibitem{ss}
A.~Schwartz and T.~Shemanske.
\newblock Maximal orders in central simple algebras and {B}ruhat-{T}its
  buildings.
\newblock {\em Journal of Number Theory}, 56(1):115 -- 138, 1996.

\bibitem{shemanske}
T.~Shemanske.
\newblock The arithmetic and combinatorics of buildings for $\text{Sp}_n$.
\newblock {\em Transactions of the American Mathematical Society}, 359(7):3409
  -- 3423, 2007.

\bibitem{woess}
W.~Woess.
\newblock {\em Random Walks on Infinite Graphs and Groups}, volume 138 of {\em
  Cambridge Tracts in Mathematics}.
\newblock Cambridge University Press, Cambridge, 2000.

\end{thebibliography}
\vfill
College of Mount Saint Vincent, 6301 Riverdale Ave., Riverdale, NY 10471\\
E-mail address: setyadi@member.ams.org
\end{document}